\documentclass[11pt]{article}
\usepackage[russian]{babel}
\usepackage[cp1251]{inputenc}
\usepackage{amsmath,latexsym}
\usepackage[psamsfonts]{amssymb}
\usepackage[dvips]{graphicx}
\usepackage[mathcal]{euscript}
\begin{document}
УДК. 517.954, 517.927.25,517.589
\begin{center}
\textbf{ Об одной задаче с условиями сопряжения для уравнения\\ четного порядка с дробной производной в смысле Капуто
 }\\
\textbf{Б.Ю.Иргашев}\\
Наманганский инженено-строительный институт.г.Наманган.Узбекистан.\\
Институт Математики Республики Узбекистан.\\
bahromirgasev@gmail.com
\end{center}

\textbf{Аннотация.} В работе для уравнения высокого четного порядка с  дробной производной в смысле Капуто изучена задача, в прямоугольной области, с условиями сопряжения. Дан критерий единственности решения. Решение построено в виде ряда Фурье по собственным функциям одномерной задачи.

 \textbf{Ключевые слова.} Уравнение четного порядка, дробная производная Капуто, разрывный коэффициент, условия сопряжения, собственное значение, собственная функция, ряд Фурье.

\textbf{1.Введение и постановка задачи.} В области $\Omega  = {\Omega _x} \times {\Omega _y},$ ${\Omega _x} = \left\{ {x:0 < x < \pi } \right\},$ ${\Omega _y} = \left\{ {y:\,  -a < y < b} \right\}$, $a ,b > 0,\,a,b \in R$, рассмотрим уравнение в частных производных
\[L\left[ u \right] \equiv \left\{ \begin{gathered}
  l\left( {u\left( {x,y} \right)} \right) + {}_CD_{0y}^\alpha u\left( {x,y} \right) = 0,\,y > 0,0 < \alpha  < 1, \hfill \\
  l\left( {u\left( {x,y} \right)} \right) + {}_CD_{0y}^\beta u\left( {x,y} \right) = 0,y < 0,1 < \beta  < 2, \hfill \\
\end{gathered}  \right.\eqno(1)\]
где
\[l\left( {u\left( {x,y} \right)} \right) = {\left( { - 1} \right)^s}\frac{{{\partial ^{2s}}u\left( {x,y} \right)}}{{\partial {x^{2s}}}} + {\left( {{p_{s - 1}}\left( x \right)\frac{{{\partial ^{s - 1}}u\left( {x,y} \right)}}{{\partial {x^{s - 1}}}}} \right)^{\left( {s - 1} \right)}} + ...\]
\[ + {\left( {{p_1}\left( x \right)\frac{{\partial u\left( {x,y} \right)}}{{\partial x}}} \right)^\prime } + {p_0}\left( x \right)u\left( {x,y} \right),\]
\[{p_j}\left( x \right) \in {C^\infty }\left( {\overline {{\Omega _x}} } \right),\,\,j = 0,1,...,s - 1,\]
$$p_j^{\left( {2i + 1} \right)}\left( 0 \right) = p_j^{\left( {2i + 1} \right)}\left( \pi  \right) = 0,\,j ,i =0,1,2... .$$
\[{}_CD_{0y}^\alpha u\left( {x,y} \right) = \frac{1}{{\Gamma \left( {1 - \alpha } \right)}}\int\limits_0^y {\frac{{\frac{{\partial u\left( {x,z} \right)}}{{\partial z}}}}{{{{\left( {y - z} \right)}^\alpha }}}} dz,{}_CD_{0y}^\beta u\left( {x,y} \right) = \frac{1}{{\Gamma \left( {2 - \beta } \right)}}\int\limits_0^y {\frac{{\frac{{{\partial ^2}u\left( {x,z} \right)}}{{\partial {z^2}}}}}{{{{\left| {y - z} \right|}^{\beta  - 1}}}}} dz\]
- дробные производные в смысле Капуто.\\
Пусть ${\Omega _ + } = \Omega  \cap \left( {y > 0} \right),{\mkern 1mu} {\mkern 1mu} {\Omega _ - } = \Omega  \cap \left( {y < 0} \right).$  Для уравнения  (1) рассмотрим следующую  задачу с условиями сопряжения.

\textbf{Задача D.}  Найти функцию $u(x,y)$  с условиями
\[Lu\left( {x,y} \right) \equiv 0,{\mkern 1mu} {\mkern 1mu} {\mkern 1mu} {\mkern 1mu} \left( {x,y} \right) \in {\Omega _ + } \cup {\Omega _ - },\]
$$\frac{{{\partial ^{2s}}u}}{{\partial {x^{2s}}}} \in C\left( \Omega  \right),\frac{{{\partial ^{2s - 1}}u}}{{\partial {x^{2s - 1}}}} \in C\left( {\bar \Omega } \right),{\mkern 1mu} {\kern 1pt}  $$
$${}_CD_{0y}^\alpha u \in C\left( {{\Omega _ + }} \right),{}_CD_{0y}^\beta u \in C\left( {{\Omega _ - }} \right),$$
$$\frac{{{\partial ^{2j}}u}}{{\partial {x^{2j}}}}\left( {0,y} \right) = \frac{{{\partial ^{2j}}u}}{{\partial {x^{2j}}}}\left( {\pi ,y} \right) = 0,{\mkern 1mu} {\kern 1pt} j = 0,1,...,s - 1; - a \leqslant y \leqslant b,\eqno(2)$$
$$u\left( {x, + 0} \right) = u\left( {x, - 0} \right),\,0 \leqslant x \leqslant \pi ,\eqno(3)$$
$$_CD_{0y}^\alpha u\left( {x, + 0} \right) = {u_y}\left( {x, - 0} \right),0 \leqslant x \leqslant \pi,\eqno(4)$$
\[u\left( {x,b} \right) - u\left( {x,a} \right) = \varphi \left( x \right),\,0 \leqslant x \leqslant \pi,\eqno(5)\]
где $\varphi \left( x \right)$ достаточно гладкая функция.

Дифференциальные уравнения в частных производных дробного порядка лежат в основе математического моделирования различных физических процессов и явлений окружающей среды, имеющих фрактальную природу [1], [2].

Уравнение (1) при $\alpha =1$, $\beta =2$, $s =1$ является параболо-гиперболическим и исследовалось во многих работах, например [3]-[6] и др.

Задачи с условиями сопряжения для уравнений с дробными производными изучались в работах [7]-[11] для уравнений второго порядка, а в работах [12],[13] для уравнений четвертого порядка c нелокальными условиями.

В данной работе изучается краевая задача с условиями сопряжения для уравнения высокого четного порядка. Найдены условия единственности решения поставленной задачи. Решение построено в виде ряда Фурье. Получены достаточные условия возможности почленного дифференцирования ряда.

\textbf{2. Единственность решения.} Учитывая (2), прежде рассмотрим следующую задачу на собственные значения:
\[\left\{ \begin{gathered}
  l\left( {X\left( x \right)} \right) = \lambda X\left( x \right), \hfill \\
  {X^{\left( {2j} \right)}}\left( 0 \right) = {X^{\left( {2j} \right)}}\left( \pi  \right) = 0,\,j = 0,1,...,s - 1. \hfill \\
\end{gathered}  \right.\eqno(6)\]
Опишем некоторые свойства собственных значений и собственных функций задачи (6). Задача (6) является самосопряженной, поэтому она имеет не более чем счетное число собственных значений ${\lambda _k},\,k = 1,2,...,$  и ортонормированную  систему собственных функций  ${X_k}\left( x \right)$ . Далее будем считать, что
$${\lambda _k} > 0,\,k = 1,2,...\,.\eqno(7)$$
условие (7) будем выполняться, если например
\[l\left( u \right) = {\left( { - 1} \right)^s}\frac{{{\partial ^{2s}}u\left( {x,y} \right)}}{{\partial {x^{2s}}}} + {p_0}\left( x \right)u\left( {x,y} \right),\,\,\,{p_0}\left( x \right) \geqslant 0.\]
Из условия (7) следует, что у задачи (6) существует функция Грина $G\left( {x,\xi } \right) = G\left( {\xi ,x} \right),$  с помощью которой задачу (6) можно свести к интегральному уравнению с симметричным ядром
\[{X_k}\left( x \right) = {\lambda _k}\int\limits_0^\pi  {G\left( {x,\xi } \right){X_k}\left( \xi  \right)d\xi } ,\]
отсюда используя неравенство Бесселя имеем
\[\sum\limits_{k = 1}^\infty  {\frac{{X_k^2\left( x \right)}}{{\lambda _k^2}} \leqslant } \int\limits_0^\pi  {{G^2}\left( {x,\xi } \right)d\xi }  < \infty .\eqno(8)\]
Известна (см. [14]) асимптотика собственных значений задачи (6):
\[{\lambda _k} = {k^{2s}} + {c_{ - 2s + 2}}{k^{2s - 2}} + ... + {c_0} + \frac{{{c_2}}}{{{k^2}}} + \frac{{{c_4}}}{{{k^4}}} + ...,\]
где постоянные ${c_j}$  вычисляются определенным образом.

Пусть теперь $u(x,y)$ некоторое решение задачи $D.$ Рассмотрим его коэффициенты Фурье по системе собственных функций $X_{k}$
\[{u_k}\left( y \right) = \int\limits_0^\pi  {u\left( {x,y} \right){X_k}\left( x \right)dx,} \]
введем функцию
\[{u_{k\varepsilon }}\left( y \right) = \int\limits_\varepsilon ^{\pi  - \varepsilon } {u\left( {x,y} \right){X_k}\left( x \right)dx,} \]
где $\epsilon >0$ -достаточно малое число, отсюда
\[{}_CD_{0y}^\alpha {u_{k\varepsilon }}\left( y \right) = \int\limits_\varepsilon ^{\pi  - \varepsilon } {{}_CD_{0y}^\alpha u\left( {x,y} \right){X_k}\left( x \right)dx}  =  - \int\limits_\varepsilon ^{\pi  - \varepsilon } {l\left( {u\left( {x,y} \right)} \right){X_k}\left( x \right)dx} .\]
Далее интегрируя по частям и переходя к пределу при $\varepsilon  \to  + 0$ , имеем
\[{}_CD_{0y}^\alpha {u_k}\left( y \right) =  - \int\limits_0^\pi  {u\left( {x,y} \right)l\left( {{X_k}\left( x \right)} \right)dx}  = \]
\[ =  - {\lambda _k}\int\limits_0^\pi  {u\left( {x,y} \right){X_k}\left( x \right)dx}  =  - {\lambda _k}{u_k}\left( y \right),\,(y > 0).\]
Аналогично будем иметь
\[{}_CD_{0y}^\beta {u_k}\left( y \right) =  - {\lambda _k}{u_k}\left( y \right),\,(y < 0).\]
Учитывая условия (3)-(5) , получим задачу
\[\left\{ \begin{gathered}
  {}_CD_{0y}^\alpha {u_k}\left( y \right) =  - {\lambda _k}{u_k}\left( y \right),\,y > 0, \hfill \\
  {}_CD_{0y}^\beta u\left( {x,y} \right) =  - {\lambda _k}{u_k}\left( y \right),y < 0, \hfill \\
  {u_k}\left( { + 0} \right) = {u_k}\left( { - 0} \right), \hfill \\
  {}_CD_{0y}^\alpha {u_k}\left( { + 0} \right) = {{u'}_k}\left( { - 0} \right), \hfill \\
  {u_k}\left( b \right) - {u_k}\left( a \right) = {\varphi _k} = \int\limits_0^\pi  {\varphi \left( x \right){X_k}\left( x \right)dx} . \hfill \\
\end{gathered}  \right.\eqno(9)\]
Решение уравнений входящих в задачу (9) будем искать в виде
\[{u_k}\left( y \right) = \left\{ \begin{gathered}
  \sum\limits_{n = 0}^\infty  {{d_n}{y^{{\gamma _1}n + {\delta _1}}}} ,\,y > 0, \hfill \\
  \sum\limits_{n = 0}^\infty  {{d_n}{{\left( { - y} \right)}^{{\gamma _2}n + {\delta _2}}}} ,\,y < 0, \hfill \\
\end{gathered}  \right.\eqno(10)\]
где ${\gamma _1},{\gamma _2},{\delta _1},{\delta _2}$ неизвестные постоянные подлежащие определению. Подставляя (10) в уравнение (9) и приравнивая коэффициенты при одинаковых степенях $y$, находим общее решение уравнения (9) в виде
\[{u_k}\left( y \right) = {c_1}{E_{\alpha ,1}}\left( { - {\lambda _k}{y^\alpha }} \right),\,(y > 0),\]
\[{u_k}\left( y \right) = {c_2}{E_{\beta ,1}}\left( { - {\lambda _k}{{\left( { - y} \right)}^\beta }} \right) + {c_3}\left( { - y} \right){E_{\beta ,2}}\left( { - {\lambda _k}{{\left( { - y} \right)}^\beta }} \right),\,(y < 0),\]
где
\[{c_1},{c_2},{c_3} - const,\]
\[{E_{\mu ,\eta }}\left( z \right) = \sum\limits_{n = 0}^\infty  {\frac{{{z^n}}}{{\Gamma \left( {\mu n + \eta } \right)}}}\eqno(11)\]
-	функция Миттаг-Леффлера.\\
Удовлетворив условиям задачи (9) имеем
\[{c_1} = {c_2},\,{\lambda _k}{c_1} = {c_3},\]
\[{c_1} = \frac{{{\varphi _k}}}{{\Delta \left( k \right)}},\]
где
\[\Delta \left( k \right) = {E_{\alpha ,1}}\left( { - {\lambda _k}{b^\alpha }} \right) - \left( {{E_{\beta ,1}}\left( { - {\lambda _k}{a^\beta }} \right) + a{\lambda _k}{E_{\beta ,2}}\left( { - {\lambda _k}{a^\beta }} \right)} \right).\]
Отсюда решение задачи (9) имеет вид
\[{u_k}\left( y \right) = \frac{{{\varphi _k}}}{{\Delta \left( k \right)}}{E_{\alpha ,1}}\left( { - {\lambda _k}{y^\alpha }} \right),\,(y > 0),\]
\[{u_k}\left( y \right) = \frac{{{\varphi _k}}}{{\Delta \left( k \right)}}{E_{\beta ,1}}\left( { - {\lambda _k}{{\left( { - y} \right)}^\beta }} \right) + \frac{{{\lambda _k}{\varphi _k}}}{{\Delta \left( k \right)}}\left( { - y} \right){E_{\beta ,2}}\left( { - {\lambda _k}{{\left( { - y} \right)}^\beta }} \right),\,(y < 0).\]
Справедлива следующая теорема единственности.

\textbf{Теорема 1.} Если существует решение задачи $D$, то оно единственно тогда и только тогда, когда выполняется условие
\[\Delta \left( k \right) = {E_{\alpha ,1}}\left( { - {\lambda _k}{b^\alpha }} \right) - \left( {{E_{\beta ,1}}\left( { - {\lambda _k}{a^\beta }} \right) + a{\lambda _k}{E_{\beta ,2}}\left( { - {\lambda _k}{a^\beta }} \right)} \right) \ne 0,\,\forall k.\]

\textbf{Доказательство.} Верность теоремы 1 следует из полноты системы собственных функций задачи (6) в пространстве ${L_2}$ [15].

Теперь нам надо показать отделимость от нуля выражения $\Delta \left( k \right)$ . Известно [16] , что функция (11) на вещественной оси, если
\[0 < \mu  < 2,\]
может иметь только конечное число нулей. В частности, если
\[1 < \mu  < 2,\,\eta  \geqslant \frac{3}{2}\mu ,\]
функция (11) не имеет вещественных нулей [17],[18].\\
Справедлива следующая лемма.

\textbf{Лемма 1.} Пусть
$$h = \max \left\{ {t:\,{E_{\beta ,2}}\left( { - t} \right) = 0,t > 0} \right\},$$
тогда существует номер $K > {k_0},$ где ${\lambda _{{k_0}}}{a^\beta } > h,$ что для всех $k > K$
выполняется оценка
\[\left| {{E_{\alpha ,1}}\left( { - {\lambda _k}{b^\alpha }} \right) - \left( {{E_{\beta ,1}}\left( { - {\lambda _k}{a^\beta }} \right) + a{\lambda _k}{E_{\beta ,2}}\left( { - {\lambda _k}{a^\beta }} \right)} \right)} \right| \geqslant \delta  > 0,\]
где $\delta $  некоторое положительное число.

\textbf{Доказательство. } Используя асимптотические разложения из [16] , при ${\lambda _k} \to  + \infty $,  имеем
\[{E_{\alpha ,1}}\left( { - {\lambda _k}{b^\alpha }} \right) = \frac{1}{{{\lambda _k}{b^\alpha }\Gamma \left( {1 - \alpha } \right)}} + O\left( {\frac{1}{{\lambda _k^2}}} \right),\]
\[{E_{\beta ,1}}\left( { - {\lambda _k}{a^\beta }} \right) = \frac{1}{{{\lambda _k}{a^\beta }\Gamma \left( {1 - \beta } \right)}} + O\left( {\frac{1}{{\lambda _k^2}}} \right),\]
\[{E_{\beta ,2}}\left( { - {\lambda _k}{a^\beta }} \right) = \frac{1}{{{\lambda _k}{a^\beta }\Gamma \left( {2 - \beta } \right)}} + O\left( {\frac{1}{{\lambda _k^2}}} \right).\]
Учитывая это получим
\[\mathop {\lim }\limits_{k \to  + \infty } \Delta \left( k \right) =  - \frac{1}{{{a^{\beta  - 1}}\Gamma \left( {2 - \beta } \right)}}.\]
Для завершения доказательства в качестве $\delta $  достаточно взять следующее число:
\[\delta  = \frac{1}{{{a^{\beta  - 1}}\Gamma \left( {2 - \beta } \right)}} - \varepsilon  > 0,\]
где  $0 < \varepsilon $ - произвольное малое число. \textbf{Лемма 1 доказана.}

Перейдем к нахождению решения задачи $D$. Учитывая вышеизложенное, формальное решение задачи $D$ имеет вид
\[u\left( {x,y} \right) = \left\{ \begin{gathered}
  \sum\limits_{k = 0}^\infty  {{X_k}\left( x \right)\frac{{{\varphi _k}}}{{\Delta \left( k \right)}}{E_{\alpha ,1}}\left( { - {\lambda _k}{y^\alpha }} \right)} ,\,y > 0, \hfill \\
  \sum\limits_{k = 0}^\infty  {{X_k}\left( x \right)\left\{ {\frac{{{\varphi _k}}}{{\Delta \left( k \right)}}{E_{\beta ,1}}\left( { - {\lambda _k}{{\left( { - y} \right)}^\beta }} \right)} \right. + }  \hfill \\
  \left. { + \frac{{{\lambda _k}{\varphi _k}}}{{\Delta \left( k \right)}}\left( { - y} \right){E_{\beta ,2}}\left( { - {\lambda _k}{{\left( { - y} \right)}^\beta }} \right)} \right\},\,y < 0. \hfill \\
\end{gathered}  \right.\eqno(12)\]
Покажем , что (12) является классическим решением поставленной задачи. Справедлива теорема.

\textbf{Теорема 2.} Пусть  выполняются условия

1. $\Delta \left( k \right) \ne 0,\,\forall k;$

2. $\varphi \left( x \right) \in {C^{4s}}\left( {{{\overline \Omega  }_x}} \right),{\varphi ^{\left( {2j} \right)}}\left( 0 \right) = {\varphi ^{\left( {2j} \right)}}\left( \pi  \right) = {\overline \varphi  ^{\left( {2j} \right)}}\left( 0 \right) = {\overline \varphi  ^{\left( {2j} \right)}}\left( \pi  \right) = 0,$

$\overline \varphi  \left( x \right) = l\left( {\varphi \left( x \right)} \right),\,j = 0,1,...,s - 1.$

Тогда (12) является классическим решением задачи $A$.

\textbf{Доказательство. } Нужно показать возможность почленного дифференцирования ряда  (12) по переменным $x,y$  до порядков входящих в уравнение (1). Формально имеем
\[l\left( {u\left( {x,y} \right)} \right) = \left\{ \begin{gathered}
  \sum\limits_{k = 0}^\infty  {{\lambda _k}{X_k}\left( x \right)\frac{{{\varphi _k}}}{{\Delta \left( k \right)}}{E_{\alpha ,1}}\left( { - {\lambda _k}{y^\alpha }} \right)} ,\,y > 0, \hfill \\
  \sum\limits_{k = 0}^\infty  {{\lambda _k}{X_k}\left( x \right)\left\{ {\frac{{{\varphi _k}}}{{\Delta \left( k \right)}}{E_{\beta ,1}}\left( { - {\lambda _k}{{\left( { - y} \right)}^\beta }} \right)} \right. + }  \hfill \\
  \left. { + \frac{{{\lambda _k}{\varphi _k}}}{{\Delta \left( k \right)}}\left( { - y} \right){E_{\beta ,2}}\left( { - {\lambda _k}{{\left( { - y} \right)}^\beta }} \right)} \right\},\,y < 0. \hfill \\
\end{gathered}  \right.\]
Далее учитывая оценку [16]
\[\left| {{E_{\mu ,\eta }}\left( { - z} \right)} \right| \leqslant \frac{M}{{1 + \left| z \right|}},\,0 < \mu  < 2,M - const,\]
при $y>0$  получим
\[\left| {l\left( {u\left( {x,y} \right)} \right)} \right| \leqslant M\sum\limits_{k = 0}^\infty  {{\lambda _k}\left| {{X_k}} \right|} \left| {{\varphi _k}} \right| \leqslant M\sqrt {\sum\limits_{k = 0}^\infty  {{{\left( {\frac{{{X_k}}}{{{\lambda _k}}}} \right)}^2}} } \sqrt {\sum\limits_{k = 0}^\infty  {{{\left( {\lambda _k^2{\varphi _k}} \right)}^2}} } .\eqno(13)\]
Первый множитель сходится за счет оценки (8). Изучим второй множитель. Имеем
\[{\varphi _k} = \int\limits_0^\pi  {\varphi \left( x \right){X_k}\left( x \right)dx}  = \frac{1}{{{\lambda _k}}}\int\limits_0^\pi  {\varphi \left( x \right)l\left( {{X_k}} \right)dx}  = \]
\[ = \frac{1}{{{\lambda _k}}}\int\limits_0^\pi  {l\left( {\varphi \left( x \right)} \right){X_k}dx}  = \frac{1}{{\lambda _k^2}}\int\limits_0^\pi  {{l^2}\left( {\varphi \left( x \right)} \right){X_k}dx} ,\]
применим неравенство Бесселя
\[\sum\limits_{k = 0}^\infty  {{{\left( {\lambda _k^2{\varphi _k}} \right)}^2}}  \leqslant \frac{1}{{\lambda _k^2}}\int\limits_0^\pi  {{{\left\{ {{l^2}\left( {\varphi \left( x \right)} \right)} \right\}}^2}dx}  < \infty .\eqno(14)\]
Из (8) и (14) следует  абсолютная и равномерная сходимость ряда (13).Аналогично доказывается сходимость ряда при $y<0$. Итак выражение (12) является классическим решением задачи $D$. \textbf{Теорема 2 доказана}.

\textbf{Замечание 1.} Из выполнения условий теоремы 2 и из теоремы Гильберта-Шмидта
 следует возможность разложения функции  $\varphi ( x )$ в ряд по собственным функциям $X_{k}$.

\textbf{Замечание 2.} Если ${\Delta }\left( k \right) = 0,$ при некоторых значениях $k = {k_1},{k_2},...{k_p}$, то для
разрешимости задачи $D$ достаточно выполнения условий ортогональности
 $\left( {{\varphi}\left( x \right),{X_k}\left( x \right)} \right) = 0,\,k = {k_1},...,{k_p}.$

\begin{center}
Литература
\end{center}
1.  А. М.Нахушев . Дробное исчисление и его применение. М., Физматлит,2003, с. 272.\\
2.  А. Н.Боголюбов , А. А. Кобликов , Д. Д. Смирнова , Н. Е. Шапкина.  Математическое моделирование сред с временной дисперсией при помощи дробного дифференцирования. Матем. моделирование, 25: 12,  2013, с. 50–64.\\
3. К. Б. Сабитов, Начально-граничная и обратные задачи для неоднородного уравнения смешанного парабологиперболического уравнения, Матем. заметки, 2017,том 102, выпуск 3, 415–435\\
4. К. Б. Сабитов, Нелокальная задача для уравнения параболо-гиперболического типа в прямоугольной области, Матем. заметки, 2011, том 89, выпуск 4, 596–602.\\
5. К. Б. Сабитов, Э. М. Сафин, Обратная задача для уравнения смешанного параболо-гиперболического типа, Матем. заметки, 2010, том 87, выпуск 6, 907–918.\\
6. К. Б. Сабитов, Задача Трикоми для уравнения смешанного параболо-гиперболического в прямоугольной области, Матем. заметки, 2009, том 86, выпуск 2, 273–279\\
7. A.S. Berdyshev, A. Cabada, E.T. Karimov. On a non-local boundary problem for a parabolic–hyperbolic equation involving a Riemann–Liouville fractional differential operator. Nonlinear Analysis: Theory, Methods and Applications. 75 (6), 2012, p. 3268-3273.\\
8. P.Agarwal, A.Berdyshev, E.Karimov. Solvability of a non-local problem with integral transmitting condition for mixed type equation with Caputo fractional derivative. Results in Mathematics. 71(3), 2017, p.1235-1257.\\
9. A.S. Berdyshev, E.T Karimov, N. Akhtaeva. Boundary value problems with integral gluing conditions for fractional-order mixed-type equation. International Journal of Differential Equations. 2011.\\
10. O. Kh. Abdullaev , K. Sadarangani  Non-local problems with integral gluing condition for loaded mixed
type equations involving the Caputo fractional derivative. Electron. J. Differential Equations, 2016.\\
11. O. Kh. Masaeva. Uniqueness of solutions to Dirichlet problems for generalized Lavrent’ev-Bitsadze
equations with a fractional derivative, Electron. J. Differential Eq., 2017 (74): 1–8.\\
12. Berdyshev A. S., Cabada A., Kadirkulov B. J. The Samarskii–Ionkin type problem for the fourth order
parabolic equation with fractional differential operator. Comput. Math. Appl., 2011. Vol. 62. P. 3884–
3893. DOI: 10.1016/j.camwa.2011.09.038 Vol. 2016. No. 164. P. 1–10. URL: https://ejde.math.txstate.edu\\
13. Berdyshev A. S., Kadirkulov J. B. On a nonlocal problem for a fourth-order parabolic equation
with the fractional Dzhrbashyan–Nersesyan operator. Differ. Equ., 2016. Vol. 52. No. 1. P. 122–127.
DOI: 10.1134/S0012266116010109\\
14. В. А. Садовничий, “О следах обыкновенных дифференциальных операторов высших порядков”, Матем. сб., 72(114):2 (1967), 293–317.\\
15. М.А. Наймарк, Линейные дифференциальные операторы, Москва,Наука,1969.\\
16. Джрбашян М. М.  Интегральные преобразования и представления функций в комплексной
области. M.: Наука.1966, с. 672.\\
17. Псху А. В. О вещественных нулях функции типа Миттаг – Леффлера. Матем. заметки, 77 (4),2005, с.592–599.\\
18. Псху А. В. Уравнения в частных производных дробного порядка. М., Наука, 2005, с. 199.

\end{document}